\documentclass[leqno,12pt]{article}[1999/09/05]
\usepackage{amsfonts,a4}
\usepackage{graphics}
\usepackage{color}
\usepackage{colordvi,epsfig}

\newtheorem{satz}{Proposition}

\newtheorem{theorem}[satz]{Theorem}

\newtheorem{corollary}[satz]{Corollary}

\newcommand{\F}[1]{ \ensuremath{ \mathbb{F}_{2^{#1}}}}
\newcommand{\Fm}[1]{ \ensuremath{ \mathbb{F}_{2^{#1}}^{\, *}}}
\newcommand{\C}{ \ensuremath{ \mathbb C}}

\newcommand{\N}{ \ensuremath{ \mathbb N}}

\newcommand{\proved}{{\hfill \bf \fbox{$ $}}}

\begin{document}
\bibliographystyle{siam}

\title{ A characterization of a class of maximum nonlinear functions}
\author{
 Doreen Hertel and Alexander Pott \\ Institute
for Algebra and Geometry\\ Otto-von-Guericke-University Magdeburg\\
D-39016 Magdeburg}

\maketitle

\begin{abstract}

Maximum nonlinear functions $F: \F{m}\to \F{m}$ are widely used 
in cryptography because the coordinate  functions 
$F_\beta(x):=tr(\beta F(x))$, $\beta \in \F{m}^\ast$, have large 
distance to 
linear functions. More precisely, the Hamming distance to
the characteristic functions of hyperplanes is large.
One class of maximum nonlinear functions are the 
 Gold power functions
$x^{2^k+1}$, $\gcd(k,m)=1$. We characterize these functions
in terms of the distance of their coordinate functions to characteristic functions of
subspaces of codimension $2$ in $\F{m}$.

\end{abstract}

\section{Introduction}
The finite field with $2^m$ elements is denoted by
$\F{m}$. The multiplicative group of the field is
denoted by $\F{m}^\ast$. We 
may also view $\F{m}$ as an
$m$-dimensional vector space over $\F{}$.
The trace function is the
linear mapping $tr:\F{m}\to\F{}$ defined
by $tr(x)=\sum_{i=0}^{m-1} x^{2^i}$. It is well known
that the mappings $tr_\beta$ defined by
$tr_\beta(x)=tr(\beta x)$ are linear, again, and all $2^m$
linear mappings $\F{m}\to\F{}$ can be represented like this.
The reader is referred to \cite{LN} for background from 
the theory of finite fields.

The {\bf  Hamming distance} between 
two {\bf boolean} functions $f, g:\F{m}\to\F{}$ is the 
number of elements $x\in\F{m}$ such that $f(x)\neq g(x)$.
The distance $d$ between $f$ and the 
linear function
$tr_\gamma$ is $2^{m-1}-\frac{t_\gamma}2$, where 
\begin{equation}\label{nl_1}
t_\gamma =\sum_{x\in\F{m}} (-1)^{tr(\gamma x + f(x))}.
\end{equation}
This is easily seen since $t_\gamma=2^m-2d$.
The distance $d$ to the affine function $x\mapsto tr_\gamma(x)+1$ 
is $2^{m-1}+\frac{t_\gamma}2$ (now $-t_\gamma=2^m-2d$). 
We say that $f$ is highly nonlinear if the smallest distance to 
all (affine) linear functions is very high. In other words, the
maximum value of $t_\gamma$ (where $\gamma\in\F{m}$) is small. 
We denote the maximum value for $|t_\gamma|$ the {\bf linearity} of $f$:
\begin{equation}\label{N(f)}
{\cal L}(f) :=\max_{\gamma \in \F{m}} 
\left| \sum_{x\in\F{m}} (-1)^{tr(\gamma x + f(x))}\right|.
\end{equation}
It is well known that 
\begin{equation}\label{lowerbound}
{\cal L}(f)\geq 2^{\frac{m}{2}}
\end{equation}
where equality may occur (of course, only for $m$ even), see
\cite{chabaud}, for instance. A function $f$ satisfying
${\cal L}(f)=2^{\frac{m}{2}}$ is called a {\bf bent function}. 
They exist for every even $m$, see \cite{chabaud}. Bent
functions are basically the same objects as certain 
difference sets, see \cite{lander-83} and \cite{bjl}.  
Note that ${\cal L}(f)=2^m$ if $f$ is linear or
affine. 

Many authors formulate this problem in a slightly different
(though equivalent) way: They define the 
{\bf nonlinearity} of $f$ to be 
\begin{equation}\label{LN}
{\cal N}(f):= 2^{m-1}-\frac12 {\cal L}(f).
\end{equation}
This number is the smallest Hamming distance between $f$ and the set
of 
all (affine) linear functions $tr(\gamma x)$ and $tr(\gamma
x)+1$ (which is called the 1st order Reed-Muller code).
The number 
${\cal N}(f)$ is the so called covering radius of this code.

The goal
is to maximize ${\cal N}(f)$. This maximum is known precisely
only in the case $m$ even, hence the 
covering radius problem is solved if $m$ is even. If $m$ is 
odd, the problem seems to be much harder.

The numbers $t_\gamma$ are sometimes called the
{\bf Walsh coefficients} ${\cal W}_f(\gamma)$. 
If $f(x)=x^d$, we simply write ${\cal W}_d(\gamma)$.
In this paper, we do not need the 
full power of the Walsh transformation, therefore we 
just refer to the literature for more information about
this important concept, see \cite{helleseth-kumar-handbook-98}, for
instance.
However, we use the notation ${\cal W}_d(\gamma)$ to denote
(\ref{nl_1})
for the power mapping $x^d$. Moreover we mention that 
the function $f$ is uniquely determined by its Walsh spectrum
$$
\{{\cal W}_f(\gamma)\, :\, \gamma\in\F{m}\}.
$$
This also holds if $f:\F{m}\to\C$ is a complex-valued function,
and we define
$$
{\cal W}_f(\gamma) := \sum_{x\in\F{m}}f(x)(-1)^{tr(\gamma x)}.
$$

We may interprete the integers ${\cal W}_d(\gamma)$ also in terms
of the intersection between certain sets. First of all
note that any function $f:\F{m}\to \F{}$ defines
a subset $D_f$ of $\F{m}$:
$$
D_f:=\{x\in\F{m}: f(x)=1\}.
$$
Vice versa, every subset gives rise to a mapping
$\F{m}\to \F{}$.

Specifically, we define for $d$ with $\gcd(d,2^m-1)=1$:
\begin{eqnarray*}
D_d:=\{x\in \F{m}: tr(x^d)=1\}\\
H^{0}(\alpha)=\{x\in\F{m}: tr(\alpha x)=0\}\\
H^{1}(\alpha)=\{x\in\F{m}: tr(\alpha x)=1\}.
\end{eqnarray*}

If $\alpha\neq 0$, the sets $H^0(\alpha)$ and $H^1(\alpha)$ are subspaces
of codimension $1$ in $\F{m}$ (hyperplanes), i.e. 
they have size $2^{m-1}$.
Therefore, we obtain for $\alpha\neq 0$ 
\begin{eqnarray}
{\cal W}_d(\alpha) & = & 2^m-2(|D_d\cap H^{0}(\alpha)|+
|(\F{m}\setminus D_d) \cap (\F{m}\setminus H^{0}(\alpha))|)\nonumber\\
& = & 
2^m-4|D_d\cap H^{0}(\alpha)|\label{translate}
\end{eqnarray}
and
\begin{equation}\label{translate2}
-{\cal W}_d(\alpha) =
2^m-4|D_d\cap H^1(\alpha)|.
\end{equation}
Since $|D_d|=2^{m-1}$ we have ${\cal W}_d(0)=0$.

Now we turn our attention to {\bf vectorial} functions $F:\F{m}\to \F{m}$.

Let $F:\F{m}\to \F{m}$ be
arbitrary. We consider the 
{\bf coordinate functions} $F_\beta(x):=tr(\beta\cdot F(x))$ from
$\F{m}$ to $\F{}$. The smallest nonlinearity of all  nonzero
coordinate functions of $F$ is
called, similar to the boolean case,
the {\bf nonlinearity} of $F$:
$$
{\cal N}(F)=\min_{\beta\in\Fm{m}}{\cal N}(F_\beta).
$$  
Similarly, the linearity is
\begin{eqnarray*}
{\cal L}(F) & = & \max_{\beta\in\Fm{m}}{\cal L}(F_\beta)\\
 & = & \max_{\beta, \gamma \in \F{m}, \beta \neq 0} 
\left| \sum_{x\in\F{m}} (-1)^{tr(\gamma\cdot x + \beta\cdot F(x))}\right|
\end{eqnarray*}
and the connection between theses two numbers is, like in (\ref{LN}),
\begin{equation}\label{LNV}
{\cal N}(F)=2^{m-1}-\frac12{\cal L}(F).
\end{equation}

Similar to the case of boolean functions $f$,
we have a (rather easy to prove) lower bound
$$
{\cal L}(F)\geq 2^{\frac{m+1}2}
$$
where equality may occur if $m$ is odd, see again
\cite{chabaud}, for instance.
Functions with ${\cal L}(F)= 2^{\frac{m+1}2}$ are 
called {\bf maximum nonlinear} or {\bf almost bent}. In this case, we
know the linearities of all the coordinate functions:
\begin{equation}\label{w-spectrum}
\sum_{x\in\F{m}} \left|(-1)^{tr(\gamma x + \beta F(x))}\right|
\in\{0,\pm 2^\frac{m+1}2\},
\end{equation}
see \cite{chabaud}, again. 
Also the multiplicities are known. The
following table shows how often the sum in (\ref{w-spectrum})
(for a fixed $\beta\ne 0$) takes the three different values:
\begin{equation}\label{w-mult}
\begin{array}{cc}
\mbox{value in (\ref{w-spectrum}) } &  \mbox{multiplicity}\\
\hline 0 & 2^{m-1}\\
2^\frac{m+1}2 & 2^{m-2}\pm 2^\frac{m-3}2\\
-2^\frac{m+1}2 & 2^{m-2}\mp 2^\frac{m-3}2.
\end{array}
\end{equation}
More precisely, if $tr(\beta F(0)) = 0$ then $2^\frac{m+1}{2}$ occurs
$2^{m-2}+ 2^\frac{m-3}2$ times, otherwise it occurs 
$2^{m-2}- 2^\frac{m-3}2$ times.

Let $k$ be an integer with $\gcd(k,m)=1$ Then the mappings
$\F{m}\to\F{m}$ defined by
$$
x^{2^k+1}\quad\mbox{or}\quad x^{2^{2k}-2^k+1}
$$
are maximum nonlinear. The first examples are 
called the {\bf Gold} power mappings, see \cite{gold-67},
the second class of mappings are the {\bf Kasami} power mappings,
see \cite{helleseth-kumar-handbook-98}. 

There are two more classes of maximum nonlinear functions known
({\bf Welch} and {\bf Niho} case), and
they are also power mappings $x^d$. They have been proved to be
maximum nonlinear only recently, see \cite{canteaut-welch} and \cite{hollmann-00}.

If $\gcd(k,m)\ne 1$ or if $m$ is even, the linearities
of the power mappings $x^{2^k+1}$ and $x^{2^{2k}-2^k+1}$
are also known, provided that $\gcd(d,2^m-1)=1$, see \cite{helleseth-kumar-handbook-98}.

It is interesting to know that, up to now, \underline{all}
maximum nonlinear mappings can be constructed from these four
classes. As usual, functions that can be constructed
from each other using some specified procedure are called
equivalent. In the case of functions, there are 
different ways to define {\it equivalence}.  A nice way to unify 
these is contained in \cite{carlet-charpin-zinoviev-98}. We refer the reader to 
\cite{buda-carlet-pott-wcc} where it is shown that the equivalence 
described in \cite{carlet-charpin-zinoviev-98} is indeed more general than the 
{\it classical} affine equivalence. When we say that all
known maximum nonlinear functions can be constructed from each other,
we do not mean that they are all {\it affine equivalent}, but
that they can be constructed from each other according to
Proposition 3 in \cite{carlet-charpin-zinoviev-98}, see \cite{buda-carlet-pott-wcc}.

In this paper we consider maximum nonlinear power functions $x^d$.
Perhaps, this class contains more mappings than just described.
In view of the connection between the Walsh coefficients
of $tr(\beta x^d)$ and the intersection 
between $D_d$ and $H^i(\alpha)$ in (\ref{translate}) and (\ref{translate2}), we obtain 
the following intersection numbers  between
$H^0(\alpha)$, $\alpha\in\Fm{m}$, and $D_d$:
\begin{equation}\label{quadric_intersect}
\begin{array}{cc}
\mbox{$|D_d\cap H^0(\alpha)|$} &  \mbox{multiplicity}\\
\hline 2^{m-2} & 2^{m-1}\\
2^{m-2}-2^\frac{m-3}2 & 2^{m-2}+ 2^\frac{m-3}2\\
2^{m-2}+2^\frac{m-3}2 & 2^{m-2}- 2^\frac{m-3}2.
\end{array}
\end{equation}

In this paper, we consider the intersection between $D_d$ 
and subspaces of codimension $2$. We can characterize the
Gold power mappings in terms of these intersection
sizes. 

If $F$ is a maximum nonlinear power mapping $x^d$, then
one can show that $x^d$ has to be a permutation, i.e.
$\gcd(d,2^m-1)=1$. It seems that this argument did 
not yet appear in the literature. It is actually due 
to Dobbertin, and it will appear in \cite{langevin-veron}. The
proof uses the fact that any maximum nonlinear mapping is an
almost perfect nonlinear function. We emphasize that not
all maximum nonlinear functions are bijective: Just by 
adding a suitable linear function $\alpha x$ one
gets nonbijective functions.

 In the case of
power mappings with $\gcd(d, 2^m-1)=1$, 
the computation of the (non)linearity simplifies.
In (\ref{p1}), put $\beta=\eta^d$ which is possible since
$\gcd(d,2^m-1)=1$. In (\ref{p2}), we replace $\gamma$ by $\alpha\eta$.
Finally in (\ref{p3}) we note that $\eta x$ runs trough
$\F{m}$ if $x$ does:

\begin{eqnarray}
{\cal L}(x^d) & = & \max_{\gamma, \beta \in \F{m}, \beta \neq 0} 
\left| \sum_{x\in\F{m}} (-1)^{tr(\gamma\cdot x + \beta\cdot x^d)}
\right|\nonumber\\
& = & \max_{\gamma, \eta \in \F{m}, \eta \neq 0} 
\left| \sum_{x\in\F{m}} (-1)^{tr(\gamma\cdot x + (\eta\cdot x)^d)}
\right|\label{p1}\\
& = & \max_{\alpha, \eta \in \F{m}, \eta \neq 0} 
\left| \sum_{x\in\F{m}} (-1)^{tr(\alpha\eta\cdot x + (\eta\cdot x)^d)}
\right|\label{p2}\\
& = & \max_{\alpha\in \F{m}}
 \left|\sum_{x\in\F{m}} (-1)^{tr(\alpha x + x^d)}\right|.\label{p3}
\end{eqnarray}

Note that 
\begin{equation}\label{p4}
\sum_{x\in\F{m}} (-1)^{tr(\alpha x + x^d)}=-1+
\sum_{x\in\F{m}^\ast} (-1)^{tr(\alpha x + x^d)}.
\end{equation}
This observation implies some connections between 
the linearity of power mappings and the crosscorrelation
between $m$-sequences and their decimations, as we will describe next.

Binary Sequences $a=(a_i)_{i\geq 0}$
($a_i\in\{0,1\}$) are called {\bf periodic} with period $n$
if $a_i=a_{i+n}$ for all $i$. The {\bf autocorrelation}
 of a binary sequence $a$
with period $n$ is defined by
\begin{eqnarray*}
  c_t(a) & := & \sum_{i=0}^{n-1} (-1)^{a_i+a_{i+t}}.
\end{eqnarray*}
The integer $t$ is called a {\bf phase shift} of the sequence $a$.
Since the sequence is $n$-periodic, we may compute the indices
modulo $n$. 

A sequence with $n$ odd and $c_t(a)=-1$ for all $1\leq t \leq n-1$
is called {\bf perfect}, see \cite{jungnickel-pott-99}
for more background on perfect sequences.  

Let $\zeta$ be a primitive element
of $\F{m}$. 
The sequences $a=(a_i)$  with 
 $a_i=tr(\zeta^i)$
 are called {\bf $m$-sequences}. They have period $2^m-1$ 
and they are perfect. Other classes of perfect
sequences are known. We refer the reader to the chapter on difference
sets in \cite{bjl} since perfect sequences
correspond to a certain class of cyclic difference sets,
see also \cite{jungnickel-pott-99}. 

Similarly to the autocorrelation, we define the
{\bf crosscorrelation} between two binary sequences
$a$ and $b$ of period $n$ by
\begin{eqnarray*}
  c_t(a,b) & := & \sum_{i=0}^{n-1} (-1)^{a_i+b_{i+t}}.
\end{eqnarray*}

Finally, we define the $d$-{\bf decimation} $a^{[d]}$ of an $n$-periodic
sequence $a=(a_i)$ by 
$a^{[d]}_i:=a_{id}$. Note that $a^{[d]}$ is an $m$-sequence
corresponding to the primitive element $\zeta^d$
if $a$ is the $m$-sequence defined by $a_i:=tr(\zeta^i)$ and $\gcd(d,2^m-1)=1$. 

If $\zeta$ is a primitive element
of $\F{m}$,  we may reformulate the righthand side of (\ref{p4}):
$$
\sum_{x\in\F{m}^\ast} (-1)^{tr(\alpha x + x^d)} = 
\sum_{i=0}^{2^m-2} (-1)^{tr(\alpha\zeta^i+ \zeta^{id})} = 
\sum_{i=0}^{2^m-2} (-1)^{tr(\zeta^{i+t}+ \zeta^{id})}
$$
(define $t$ by $\alpha=\zeta^t$). 
This shows that 
the linearity of the power mapping $x^d$ is the
same as $-1$ plus the maximum crosscorrelation value
between an $m$-sequence and its $d$-decimation.

\section{Main Theorem}

\begin{theorem}\label{satz1}
Let $m$ be odd and let 
$x^d$ be a maximum nonlinear power function on $\F{m}$.
Let
$$
 H^{i,j}(\alpha,\beta) := \{ x\ :\  tr(\alpha x)=i,\ 
tr(\beta x)=j \}.  
$$
Then $d=2^k+1$ for some integer $k$ with $\gcd(k,m)=1$
(i.e. $d$ is a Gold exponent) \underline{if and only if} 
\begin{equation}\label{inter2}
 |H^{i,j}(\alpha,\beta)\cap D_d |\in  \{ 2^{m-3},2^{m-3} \pm 2^{\frac{m-3}{2}} \}
\end{equation}
for all $\alpha, \beta \in \F{m}^\ast$, $\alpha\neq\beta$, and $i,j\in\F{}$.
\end{theorem}

The sets $H^{i,j}(\alpha,\beta)$,  $\alpha\neq\beta$, 
$\alpha, \beta \in \F{m}^\ast$, are precisely the 
subspaces of dimension $m-2$ in $\F{m}$. 
The set $D_d$ has some interesting properties. It is the
set of $2^{m-1}$ points in the $m$-dimensional vector space $\F{m}$
over $\F{}$. If $d$ is a Gold exponent, this set 
is a non-degenerate quadric, see \cite{helleseth-kumar-handbook-98}, 
for instance.
If $m$ is odd, there is up to equivalence only one non-degenerate 
quadric in $\F{m}$, and the intersection between this quadric
and subspaces of codimension $2$ must be the three values described
in (\ref{inter2}). This is well known to geometers, see
\cite{hirschfeld-98}, for instance. It follows from \cite{games-86} that the
only quadrics corresponding to the coordinate functions
of maximum nonlinear power mappings
are nondegenerate: In \cite{games-86}, the intersection sizes 
between quadrics $Q$ and hyperplanes are determined whenever 
$|Q|=2^{m-1}$. This applies to the situation of maximum
nonlinear power mappings $x^d$ and $m$ odd, since in this case
$|D_d|=2^{m-1}$ (because $x^d$ is bijective). It 
turns out that the intersection sizes in \ref{quadric_intersect} 
occur only in the nondegenerate case.

It is natural to ask whether there are  values
$d$ such that $D_d$ is not a nondegenerate quadric 
but has the same intersection sizes with hyperplanes. 
These
objects are called  by geometers {\bf quasi-quadrics}. Many examples
of quasi-quadrics are known, see \cite{clerck-2000}, for instance. Note that
all maximum nonlinear power mappings yield quasi-quadrics. 
Our research was motivated by the question whether 
 the quasi-quadrics constructed from maximum nonlinear functions
may also behave like quadrics if the intersection sizes with subspaces of
codimension $2$ are considered. The answer, given by
Theorem \ref{satz1}, is \underline{no}.

An interesting corollary is the following:

\begin{corollary}
The only maximum nonlinear power mappings $x^d$ on $\F{m}$
such that $D_d$ is a quadric are the Gold power mappings.
\end{corollary}

Before we are going to prove our Theorem, let us mention
the following Proposition which may be of interest in its own:

\begin{satz}\label{lemma1}
Let $x^d$ be a maximum nonlinear power mapping on $\F{m}$ with 
$gcd(d,m)=1$.
Then 
\begin{equation}
  |\ H^{i,j}(\alpha,\beta) \cap D_d\ | \in 
\{ 2^{m-3} + h\cdot 2^{\frac{m-5}{2}}\ :\ \quad -3\leq h \leq 3\},
\end{equation}
where $\alpha,\beta\in\F{m}^\ast$, $\alpha\neq \beta$.
\end{satz}

{\bf Proof.}
We define
$$
S^{i,j}(\alpha,\beta)=| H^{i,j}(\alpha,\beta) \cap D_d |
$$
and
$$
S^{i}(\alpha)=| H^i(\alpha)\cap D_d |.
$$
Assume $\alpha\neq\beta$, $\alpha,\beta\in\F{m}^\ast$. We obtain
\begin{eqnarray*}
|{\cal W}_d(\alpha+\beta)| & = & |\sum\limits_{x\in \F{m}} (-1)^{tr(\alpha
  x+\beta x+x^d)}|\\
& = & |\sum\limits_{x\in \F{m}} (-1)^{tr(\alpha
  x+i+\beta x+j+x^d)}|\\
& = & 2^m-2(3\cdot 2^{m-1}-2S^i(\alpha)-2S^j(\beta)-2\cdot
2^{m-2}+4S^{i,j}(\alpha, \beta))\\
& = & -2^m+4S^i(\alpha)+4S^j(\beta)-8S^{i,j}(\alpha, \beta).
\end{eqnarray*}
Because of (\ref{translate}) and (\ref{translate2}), we have
$$
|{\cal W}_d(\alpha+\beta)|=2^m\pm{\cal W}_d(\alpha)\pm{\cal W}_d(\beta)-
8S^{i,j}(\alpha, \beta)),
$$
hence
\begin{equation}\label{7werte}
S^{i,j}(\alpha, \beta) = 2^{m-3}\pm \frac18(\pm{\cal W}_d(\alpha+\beta)
\pm {\cal W}_d(\alpha) \pm {\cal W}_d(\beta)).
\end{equation}
This shows that there are only the seven possible values
for $S^{i,j}(\alpha, \beta)$ stated in the Proposition.\proved

The proof of Theorem \ref{satz1} reduces to the proof of an interesting
property of the trace function. This Theorem has been independently
obtained by Ph. Langevin and P. V\'{e}ron \cite{langevin-veron}. 
The proof given in their paper 
is different from ours. The Langevin-V\'{e}eron proof 
is shorter and more elegant, though less elementary. 

Theorem \ref{satz2} is not true
in the case $m$ even. For instance, if we take $m=8$ and $d=51$, then
$tr(x^d+(x+1)^d+1)=0$ for all $x\in\F{8}$. This example can be extracted
from the proof, since line (8) of the following algorithm 
does not produce the desired element $w$ if $m$ is even.
Note that $\gcd(51,2^8-1)\ne1$. If we restrict ourselves to the
case $\gcd(d,2^m-1)=1$, Theorem \ref{satz2} remains true also if
$m$ is even, see Lemma 2 in \cite{kyureghyan-wcc}.

 \begin{theorem}\label{satz2}
 Let $m$ be odd and $d \in \{ 2,...,2^m-2 \} $ odd. We have
 \begin{eqnarray*}
   tr(x^d+(x+1)^d+1) & = & 0
 \end{eqnarray*}
 for all $x \in \F{m}$, if and only if  $d=2^k+1$ for one $k \in \N $.
 \end{theorem}
 
We postpone the proof of Theorem \ref{satz2}  to the next Section. 
We are now going to show that it is sufficient to prove Theorem
\ref{satz2} in order to check Theorem \ref{satz1}.

Let $x^d$ be a maximum nonlinear power function on $\F{m}$,
hence the Walsh spectrum
$\{ {\cal W}_d(\alpha)| \alpha\in \F{m} \}$ contains only the three value 
$ \pm 2^{\frac{m+1}{2}}  $
and $0$.
We define the function $b:\F{m}\to\F{}$ by 
$$b(\alpha) = \left\{\begin{array}{cl}
1 & \mbox{if\ } {\cal W}_d(\alpha)\neq 0\\
0  & \mbox{otherwise.}
\end{array}\right.$$
If all or precisely one of the values ${\cal W}_d(\alpha)$, ${\cal W}_d(\alpha)$ and 
${\cal W}_d(\alpha)$ in equation (\ref{7werte}) are  $\neq 0$, it is impossible that
$S^{i,j}(\alpha, \beta)\in \{ 2^{m-3},2^{m-3} \pm 2^{\frac{m-3}{2}} \}.
$
Therefore,
$b(\alpha)+b(\beta)=b(\alpha+\beta)$,
hence $b$ is linear and therefore 
$$
b(x)=tr(\gamma x)\quad(=tr_\gamma(x))
$$
for some $\gamma\in\F{m}^\ast$.
If we think of $tr(x)$ as an elment in $\C$, we obtain
\begin{eqnarray}
{\cal W}({tr_\gamma})(\omega) & = & \sum_{x\in\F{m}} tr(\gamma
  x)\cdot (-1)^{tr(\omega x)}\nonumber\\
                     & = & \sum_{x \in \F{m}, tr(\gamma (x))=1} 
(-1)^{tr(\omega x)}\nonumber \\
  \label{version_1}                   & = & \left\{ \begin{array}{cl} 
                                   -2^{m-1} & \mbox{ if } \omega=\gamma \\
                                    2^{m-1} & \mbox{ if } \omega=0 \\
                                    0       & \mbox{ otherwise. }
                           \end{array} \right.
\end{eqnarray}

On the other hand, the function $b$ satisfies
$$
b(x) = \frac{1}{2^{m+1}}[{\cal W}_d(x)]^2.
$$
We compute the Walsh transform again:
\begin{eqnarray*}
  {\cal W}(b)(\omega) & = & \sum_{x\in \F{m}} \frac{1}{2^{m+1}} \left( {\cal
      W}_d(x) \right)^2 (-1)^{tr(\omega x)} \\
                 & = & \frac{1}{2^{m+1}} \sum_{x,y,z \in \F{m}} (-1)^{tr(x(y+z+\omega)+(y^d+z^d))} \\
                 & = & \frac{1}{2^{m+1}} \sum_{y,z \in \F{m}} (-1)^{tr(y^d+z^d)}
                       \hspace{-2cm}\underbrace{\sum_{x\in \F{m}} (-1)^{tr(x(y+z+\omega))}}_{
                       \footnotesize{\hspace{4cm}=\left\{ \begin{array}{cl} 2^m & \mbox{ if } z=\omega+y \\
                                   0 & \mbox{ otherwise}\end{array} \right. }} \\
                 & = & \frac{1}{2} \sum_{y\in \F{m}} (-1)^{tr(y^d+(y+\omega)^d)}  .
\end{eqnarray*}
We compare this with (\ref{version_1}) and obtain
\begin{eqnarray}\label{unique}
  \sum_{y\in \F{m}} (-1)^{tr(y^d+(y+\omega)^d)}& = & \left\{ \begin{array}{cl}
                        -2^m & \mbox{ if } \omega=\gamma \\
                        2^m & \mbox{ if } \omega=0 \\
                        0 & \mbox{ otherwise.}
                      \end{array} \right.
\end{eqnarray}
The case $\omega=\gamma $ implies
\begin{eqnarray} \label{spurwert}
  tr(y^d+(y+\gamma )^d)  = 1 && \mbox { for all } y\in \F{m} .
\end{eqnarray}
We can show that necessarily $\gamma =1$:
\begin{eqnarray*}
  tr((y+\gamma )^d) & \stackrel{(\ref{spurwert})}{=} & tr(y^d)+1 \\
                   & = & tr(y^{2^l d})+1 \\
                   & \stackrel{(\ref{spurwert})}{=} & tr((y^{2^l}+\gamma )^d) \\
                   & = & tr( (y+\gamma^{2^{m-l}})^d )
\end{eqnarray*}
for all $l=0,...,m-1$ and $y\in \F{m}$.
Suppose, that $\gamma \ne \gamma ^{2^k}$ for some $k\in \N$, then
$tr(y^d+(y+\gamma ^{2^k})^d)=tr(y^d+(y+\gamma )^d)=1$, thus
$\sum_{y\in \F{m}} (-1)^{tr(y^d+(y+\gamma ^{2^k})^d)} =\sum_{y\in \F{m}} (-1)^{tr(y^d+(y+\gamma )^d)}=-2^m$.
This is a  contradiction to the uniqueness of $\gamma$. 
Thus we have $\gamma ^{2^l}=\gamma$
for all $l=0,...,m-1$, and therefore $\gamma =1$.

Since $m$ is odd we have $tr(1)=1$. Therefore
\begin{eqnarray} \label{nullspur}
  tr(y^d+(y+1)^d+1) & = & 0
\end{eqnarray}
for all $y\in \F{m}$.
This shows that it is enough to prove Theorem \ref{satz2}.

\section{Proof of Theorem \ref{satz2}}

If $d$ satisfies
 (\ref{nullspur}), then each $d' \in D$ with $D := \{ \, 2^i d \bmod (2^m-1)\, : \, i=0,...,m-1 \, \}$
 also satisfies (\ref{nullspur}).
 We choose the smallest odd $d'$ in $D$, 
 and from now on, we denote this element by $d$.
 
 Let $w(a)$ be the binary weight of $a$.
 If $a=\sum_{i=0}^{n} z_i2^i$ is the binary representation
 of $a$, we denote the vector $(z_n,\ldots,z_0)$ by $\overline{a}$.
 We have $\sum_{i=0}^{n}z_i=w(a)$. In the following proof, all integers $a$ that occur are
 $\leq 2^{m}-1$, i.e. $\overline a$ is a vector of length at most $m$.
 By adding $0$'s, if necessary, we assume that $\overline a$ is always a vector
 of length $m$. 
 Let $v=(v_{m-1},...,v_0)$ be a vector of length $m$.
 We denote by $v^{(t)}$ the cyclic shift of the vector $v$
 about $t$ positions to the left, i.e. $v^{(t)}=(v_{m-t-1},...,v_{0},
 v_{m-1},...,v_{m-t})$.
 
 Let $d'=2^id \bmod (2^m-1)$, then $\bar{d'} = \bar{d}^{(t)}$,
 in particular $w(d')=w(d)$.
 
We define the polynomial $p$ by
\begin{eqnarray*}
  p(x) & := & x^d+(x+1)^d+1
\end{eqnarray*}
and $q$ by
\begin{eqnarray*}
q(x) & := & tr(p(x)) = \sum_{i=0}^{m-1} (p(x))^{2^i}.
\end{eqnarray*}
We have  $q(0)=0$, therefore we have to show that
\begin{equation}\label{null_q}
  q(\alpha)=0\quad\mbox{for all }\alpha\in\Fm{m}.
\end{equation}
 
Let $T=\{t_1,\ldots, t_n\}$ denote the set of exponents  which occur in $p$.
We define the set
$$
T(t)=\{\, 0\leq s\leq 2^m-2\ :\ \overline{s}^{(i)}=
\overline{t},\ i=0,\ldots,m-1\}.
$$
We obtain
\begin{eqnarray*}
q(x) & = &\sum_{t\in T} \sum_{s\in T(t)} x^s.
\end{eqnarray*}
 
In order to prove (\ref{null_q}), we have to show that
every exponent occurs an even number of times in $q(x)$.

Let $\overline{d}=(z_{m-1},\ldots,z_0)$ be the
binary vector corresponding to $d$. Since $d$ is odd,
we have $w(d)\neq 1$. If $d=2^k+1$ is a Gold exponent, then
$w(d)=2$ and $q(x)$ satisfies (\ref{null_q}) (note that
$p(x)=x^{2^k}+x$ in this case).
Hence we may assume $w(d)\geq 3$.

If $w(d)=3 $, then $d=2^k+2^l+1$ and $k>l>0$.
For the polynomials $p$ and $q$ we obtain
\begin{eqnarray*}
  p(x) & = & x^{2^k+2^l} + x^{2^k+1} + x^{2^l+1} + x^{2^k} + x^{2^l} + x \\
  q(x) & = & \sum_{i=0}^{m-1} \Big( (x^{2^k+2^l})^{2^i} + (x^{2^k+1})^{2^i} + (x^{2^l+1})^{2^i} + x^{2^i} \Big) .
\end{eqnarray*}
In $p(x)$, the exponents of binary weight $1$ (resp. $2$) occur three
times, therefore we have an odd number of
exponents of weight $1$ (resp. $2$) in $q(x)$, and therefore $q(x)$ cannot satisfy
(\ref{null_q}). This argument can be generalized: If $z=w(d)$ then
there are precisely ${z}\choose{i}$ exponents $t$ in $p(x)$
with $w(t)=i$, $1\leq t\leq d-1$.
Note that $x^d$ and 1 do not occur in $p(x)$.
If $z$ is not a power of $2$,
at least one of these binomial coefficients is odd. Therefore,
we only have to consider the case $z=2^n$, $n>1$.

Let $v$ be a binary vector of length $m$.
A {\bf subvector} $w=(w_{m-1},...,w_0)$ of $v$ is a binary vector $w\neq 0,v$ of length $m$
such that $v_i=0$ implies $w_i=0$.
  
  
  
The set of all subvectors of $\bar{d}$ is
the set of the binary vectors of the exponents
that occur in  $p(x)$.
  
In order to show that (\ref{null_q}) holds, we have to prove that
the cardinality of the multiset
\begin{eqnarray*}S(s) & := & \{ \bar{s}^{(i)}\ : \ \bar{s}^{(i)}
  \mbox{ subvector of } \bar{d},\  0\leq  i \leq m-1   \}
\end{eqnarray*}
is even for all $s \in T$. Note that it is possible that 
$\bar{s}^{(i)}=\bar{s}^{(j)}$ for $i\ne j$.

We define a {\bf gap} to be a substring $v$
of the form
$0...0$. The number $s$ of $0$'s in this substring is called the
length of the gap,
similarly for {\bf runs} which are substrings of the form $1...1$.
More precisely: There is an $i$ such that $z_{i}=z_{i+1}=\dots=z_{i+s-1}=0$, where the indices
are computed modulo $m$, i.e. we view $\bar{d}$ as a ``cyclic'' vector.
If $v=(v_iv_{i+1}\ldots v_j)$  is a substring,
we say that the indices $i,\ldots,j$ are contained in $v$.

By the following algorithm we construct a subvector $w$ of $\bar{d}$
such that $|S(w)|$ is odd. Therefore, (\ref{null_q}) is not satisfied.

{\bf Algorithm}

 Input: binary vector $\bar{d}=(z_{m-1},...,z_0)$ of weight $z=2^n,
 n\in \N$,
$m$ odd\\
 Output: subvector $w$ of $\bar{d}$ such that $|S(w)|$ is odd
 
 \begin{tabular}{rl}
 (1) & $l:= \mbox{ maximum length of a run in } \bar{d}$; \\
     & $s:= \mbox{ multiplicity of a run of length $l$ in $\bar{d}$} $;\\
     & $v:= \mbox{ run of length } l;$ \\
     & $s_{old}:=m+1$; $x_{old}:=0$;\\
 (2) & while ($w$ is not defined) do \\
 (3) & \mbox{~}\hspace{0.5cm} $y:=(y_{m-1},...,y_0)$ with \\
     & \mbox{~}\hspace{1.2cm}
       $y_i=\left\{\begin{array}{cl} 1 &
                \mbox{if $i$ is contained in a substring $v$ and $z_i$ is 1}\\ 0 & \mbox{otherwise.} \end{array} \right. $\\
 (4) & \mbox{~}\hspace{0.5cm} if $z\ne l\cdot s$ then $w:=y$; end if; \\
 (5) & \mbox{~}\hspace{0.5cm} if $z = l\cdot s$ then \\
     & \mbox{~}\hspace{1.2cm} $x :=$ minimum length of a gap between two 
 substrings $v$ in $y$; \\
     & \mbox{~}\hspace{1.2cm} $L:=$ gap of length $x$; \\
     & \mbox{~}\hspace{1.2cm} if $s=1$ then  \\
 (6) & \mbox{~}\hspace{1.9cm} if $s_{old}=m+1$ then $w:= \bar{d}-(0...010)$; end if; \\
 (7) & \mbox{~}\hspace{1.9cm} if $s_{old}\ne m+1$ then $w:=(0...0v_{old}L_{old} v_{old})$; end if; \\
     & \mbox{~}\hspace{1.2cm} end if; \\
 (8) & \mbox{~}\hspace{1.2cm} if $s=2$ then $w:=(0...0vLv)$; end if;\\
 (9) & \mbox{~}\hspace{1.2cm} if $s>2$ then \\
     & \mbox{~}\hspace{1.9cm} $s_{old}:=s;\, l_{old}:=l;\,x_{old}:=x;\,L_{old}:=L;\,v_{old}:=v;$  \\
     & \mbox{~}\hspace{1.9cm} let $v$ denote a substring of type $(v_{old}Lv_{old}...Lv_{old})$ in $\bar{d}$ \\
     & \mbox{~}\hspace{2.6cm} with maximum number $l$ of 1's; \\
     & \mbox{~}\hspace{1.9cm} $v:=(v_{old} L v_{old} ... L v_{old})$ with $v_{old}$ occurs $l/l_{old}$-times;  \\
     & \mbox{~}\hspace{1.9cm} $s:=$ multiplicity of $v$ in $\bar{d}$;\\
     & \mbox{~}\hspace{1.2cm} end if;\\
     & \mbox{~}\hspace{0.5cm} end if;\\
     & end while;\\
 \end{tabular}

 The algorithm terminates if  $z\ne l\cdot s$ or $s \leq 2$.
 Note, if the case $z\ne l\cdot s$ does not occur then  such an $s$
exists
 because $0<s<s_{old}$ in each step of the algorithm.

 Line (4): If $z\ne l\cdot s$, i.e. $y\ne \bar{d}$
 and $w=y$ is a subvector of $\bar{d}$.
 We  have $|S(w)|=1$, because
 none of the cyclic shifts $w^{(t)}\ne w$ is a subvector of $\bar{d}$.
 Suppose the vector $w^{(t)}$ with $w^{(t)} \ne w$ is a subvector of $\bar{d}$.
 Note, that $w$ and $w^{(t)}$ have the same number of 1.
 If $w^{(t)}\ne w$, then there exists a $1$ in $\bar{d}$ and this $1$
 is in $w^{(t)}$ and not in $w$.
 Because $w^{(t)}$ is a cyclic shift of $w$, this $1$ is in a string $v$, therefore
 this $1$ is in $w$. This is a contradiction to the definition of $w$.

 Line (5): If $z= l\cdot s$, then $l=2^{l'}$ and $s=2^{s'}$.
 We call the gaps $L_{j}$, $j=1,...,s$ between the runs $v$.
 Now we know, that $\bar{d}$ has the form
  \begin{eqnarray*}\bar{d}&= &(L_{s}vL_{s-1}v...L_{2}vL_{1}v).\end{eqnarray*}
 If $s > 1$, the number of gaps is even. Since $m$ is odd,
 the number of gaps with odd length and the number of gaps with even
 length is odd.
 Therefore the maximum and minimum gap have different length.
 Note by the choice of $d \in D$ to be odd it follows that 
 $L_{s}$ is one of the maximum gaps
 and has length $>x$, the mimimum length of a gap.

 Line (6): If $z= l\cdot s$ with $s=1$ and $s_{old}=m+1$ then $l \geq 4$ and $\bar{d}=(0...01...1)$.
 For $w=\bar{d}-(0...010)$ we have $|S(w)|=1$.

 Line (7): If $z= l\cdot s$ with $s=1$ and $s_{old} \ne m+1$, then $s_{old} \geq 4$.
 The vector $\bar{d}$ has the form
  \begin{eqnarray*}\bar{d} & = & (L_{s} v) \ = \ ( L_{s_{old}}  v_{old}   L   v_{old}  ...   L   v_{old}   L   v_{old}) ,\end{eqnarray*}
 where $L$ is the gap of length $x_{old}$.
 We obtain $|S(w)|=s_{old}-1$ is odd.

Line (8):
If $s=2$, then is $\bar{d}=(L_{2} v L_{1} v)$. The gap
$L_{2}$ is longer than the gap $L_{1}$.
It is easy to see that $|S(w)|=1$.

 Line (9): The new initialisation for the next while loop.
 \proved
 
 We illustrate the algorithm with an example. Here we have
 $m=23$ and $d=1+2^2+2^4+2^7+2^9+2^{11}+2^{15}+2^{17}$.
 
 Input: $\bar{d}=(0 00 00 10 10 00 10 10 10 01 01 01)$\\
 \begin{tabular}{rl}
 (1) & $z:=8$; $l:= 1$; $s:= 8$; $v:= 1$; $s_{old}:=24$; $x_{old}:=0$;\\
 (3) & $y:=\bar{d}$;  \\
 (5) & $x := 1$ ; $L:=0$; \\
 (9) & $s_{old}:=8;\,l_{old}:=1;\,x_{old}:=0;\,L_{old}:=0;\,v_{old}:=1;$ \vspace{0.1cm} \\
     & \mbox{~}\hspace{2.0cm}
        $y=(0 00 00\underline{10 1} 0 00\underline{10 10 1} 0 0\underline{1 01 01});$
       \vspace{0.1cm} \\
     & $l:=3$; $v:= 1 0 1 0 1 $;  $s:=2$; \\
 (3) & $y:=(0 00 00 00 00 00 v 0 0 v)$;\\
 (4) & $w:=y$;
 \end{tabular}\\
 Output: $w:=(0 00 00 00 00 00 10 10 10 01 01 01)$.

 \def\cprime{$'$} \def\cprime{$'$} \def\cprime{$'$}
  \def\Dbar{\leavevmode\lower.6ex\hbox to 0pt{\hskip-.23ex \accent"16\hss}D}
  \def\Dbar{\leavevmode\lower.6ex\hbox to 0pt{\hskip-.23ex \accent"16\hss}D}


\end{document}